\documentstyle[11pt,math-art]{article}

\newcommand{\cG}{{\cal G}}

\newcommand{\suchthat}{\mathchoice{\;\;|\;\;}{\;|\;}{\;|\;}{\;|\;}}

\newcommand{\mfn}[1]{{\rm #1}}

\newtheorem{theorem}{Theorem}[section]

\newtheorem{lemma}[theorem]{Lemma}
\newtheorem{corollary}[theorem]{Corollary}

\newtheorem{examplehidden}[theorem]{Example}

\newtheorem{definitionhidden}{Definition.}

\newenvironment{example}{\begin{examplehidden}\rm}%
{\end{examplehidden}\par}
\newenvironment{definition}{%
\begin{definitionhidden}\rm}%
{\end{definitionhidden}\par}
\newcommand{\proof}{\paragraph*{\hspace*{-\parindent}Proof.}}

\newcommand{\noindpar}[1]{\paragraph*{\hspace*{-\parindent}#1.}}
\newcommand{\qed}{\hspace*{\fill}\rule{2.5mm}{2.5mm}%
\vspace*{8pt}\par}
\newsavebox{\cardofbox}
\newlength{\cardofboxheight}
\newlength{\cardofboxdepth}
\newlength{\cardoftmpld}
\sbox{\cardofbox}{$l$}
\newlength{\cardoftmpls}
\sbox{\cardofbox}{$\scriptstyle l$}
\newlength{\cardoftmplss}
\sbox{\cardofbox}{$\scriptscriptstyle l$}
\newcommand{\cardof}[1]{%
\mathchoice{%
\protect\sbox{\cardofbox}{$\displaystyle{#1}\rule{0pt}{\cardoftmpld}$}%
\setlength{\cardofboxdepth}{-\dp\cardofbox}%
\addtolength{\cardofboxdepth}{-2pt}%
\setlength{\cardofboxheight}{\dp\cardofbox}%
\addtolength{\cardofboxheight}{\ht\cardofbox}%
\addtolength{\cardofboxheight}{3pt}%
\,\rule[\cardofboxdepth]%
{.5pt}{\cardofboxheight}%
\,\usebox{\cardofbox}\,%
\rule[\cardofboxdepth]{.5pt}{\cardofboxheight}\,}%
{%
\protect\sbox{\cardofbox}{${#1}\rule{0pt}{\cardoftmpld}$}%
\setlength{\cardofboxdepth}{-\dp\cardofbox}%
\addtolength{\cardofboxdepth}{-2pt}%
\setlength{\cardofboxheight}{\dp\cardofbox}%
\addtolength{\cardofboxheight}{\ht\cardofbox}%
\addtolength{\cardofboxheight}{3pt}%
\,\rule[\cardofboxdepth]%
{.5pt}{\cardofboxheight}%
\,\usebox{\cardofbox}\,%
\rule[\cardofboxdepth]{.5pt}{\cardofboxheight}\,}
{%
\protect\sbox{\cardofbox}{$\scriptstyle{#1}\rule{0pt}{\cardoftmpls}$}%
\setlength{\cardofboxdepth}{-\dp\cardofbox}%
\addtolength{\cardofboxdepth}{-1pt}%
\setlength{\cardofboxheight}{\dp\cardofbox}%
\addtolength{\cardofboxheight}{\ht\cardofbox}%
\addtolength{\cardofboxheight}{1,5pt}%
\rule[\cardofboxdepth]%
{.5pt}{\cardofboxheight}%
\hspace*{.5pt}\usebox{\cardofbox}\hspace*{.5pt}%
\rule[\cardofboxdepth]{.5pt}{\cardofboxheight}}
{%
\protect\sbox{\cardofbox}{$\scriptscriptstyle{#1}\rule{0pt}{\cardoftmplss}$}%
\setlength{\cardofboxdepth}{-\dp\cardofbox}%
\addtolength{\cardofboxdepth}{-.5pt}%
\setlength{\cardofboxheight}{\dp\cardofbox}%
\addtolength{\cardofboxheight}{\ht\cardofbox}%
\addtolength{\cardofboxheight}{.8pt}%
\rule[\cardofboxdepth]%
{.5pt}{\cardofboxheight}%
\hspace*{.5pt}\usebox{\cardofbox}\hspace*{.5pt}%
\rule[\cardofboxdepth]{.5pt}{\cardofboxheight}}%
}

\begin{document}
\title{Notes on The Connectivity of Cayley Coset Digraphs}
\author{Emanuel Knill\thanks{This work was performed under the
auspices of the U.S. Department of Energy under Contract No. W-7405-ENG-36.}\\
Los Alamos National Laboratory, CIC-3\\
Mailstop B265, Los Alamos, NM 87545\\
Email: knill@lanl.gov}
\date{January 1993}

\maketitle

\noindpar{Abstract}
Hamidoune's connectivity
results~\cite{hamidoune:hierarchical} for hierarchical Cayley digraphs are extended to
Cayley coset digraphs and thus to arbitrary
vertex transitive digraphs.
It is shown that if a Cayley coset digraph can be hierarchically
decomposed in a certain way, then it is optimally vertex connected.
The results are obtained by extending the methods used
in~\cite{hamidoune:hierarchical}.
They are used to show that cycle-prefix graphs~\cite{chen:basicCP}
are optimally vertex connected. This implies that cycle-prefix graphs
have good fault tolerance properties.

\section{Introduction}
\label{section:introduction}

Good interconnection networks for parallel computing usually have
the following properties~\cite{bermond:fault,akers:group}:
They are symmetric, so that each node
has the same view of the network. There are simple routing methods for
finding paths from one node to another. The number of edges is small.
The maximum distance between two nodes is small. The network
can be easily constructed in 2 or 3 dimensions. The network is fault
tolerant.

This motivates the study of vertex transitive, small degree
and diameter, optimally connected digraphs.
Sabidussi~\cite{sabidussi:thecomposition} showed that the class of vertex
transitive digraphs is the same as that of Cayley coset digraphs.
For that reason, there has been a substantial effort to construct
and analyze Cayley coset digraphs with good interconnection network
properties~\cite{baumslag:fault,chen:basicCP,akers:group}.

The connectivity properties of vertex transitive graphs were studied
by Watkins~\cite{watkins:connectivity} and 
Mader~\cite{mader:uberden,mader:minimale}; and of digraphs by
Hamidoune%
~\cite{hamidoune:sur,hamidoune:connectivity,hamidoune:connectivite}.
Mader~\cite{mader:minimale} showed that connected vertex transitive
digraphs have optimal edge connectivity.  The first general results on
the vertex connectivity of connected vertex transitive graphs were
obtained by Mader~\cite{mader:uberden} and
Watkins~\cite{watkins:connectivity}.  They show that every connected
edge and vertex transitive graph has optimal vertex connectivity
(i.e. the vertex connectivity is the same as the degree).
Mader~\cite{mader:uberden,mader:minimale} shows
that every connected vertex transitive graph without $K_4$ is
optimally vertex connected. He also shows that every connected edge
transitive graph is optimally vertex connected.  This work was
extended to Cayley digraphs by Hamidoune~\cite{hamidoune:sur,%
hamidoune:connectivity,hamidoune:connectivite,%
hamidoune:hierarchical}.  In~\cite{hamidoune:connectivite}, the
abelian Cayley digraphs without $K_4$ which are not optimally
vertex connected are characterized.  In~\cite{hamidoune:hierarchical}, it is
shown that connected {\em hierarchical\/} Cayley digraphs are
optimally vertex connected.  A similar
result is obtained in Baumslag~\cite{baumslag:fault} using more direct
methods.

The main tool for obtaining connectivity results in vertex transitive
graphs is the concept of an {\em atom\/}, which, briefly, is a minimal
part of the graph with connectivity many neighbors.  In this note,
atoms are used
for proving vertex connectivity results for Cayley
coset digraphs.  The main result is
Theorem~\ref{theorem:decomposition}, which generalizes Proposition~3.1
of~\cite{hamidoune:hierarchical}.  This result is used to obtain a
hierarchical decomposition result for Cayley coset digraphs, which as
a corollary yields the result of
Hamidoune~\cite{hamidoune:hierarchical} and
Baumslag~\cite{baumslag:fault} that connected hierarchical Cayley
graphs are optimally vertex connected. The main result is
applied to show that the cycle-prefix graphs which were proposed
as interconnection networks in~\cite{chen:basicCP} are optimally
vertex connected.

This note is organized as follows.  Caley coset digraphs are defined
in Section~\ref{section:caley coset}, which also contains some
elementary observations on Caley coset digraphs.  Vertex connectivity
is discussed in detail in Section~\ref{section:vertex connectivity}.
Most of this section follows closely the methods described
in~\cite{hamidoune:hierarchical} but generalizes them to Caley coset
digraphs. In Section~\ref{section:applications}, the main
result is applied to hierarchical Cayley coset digraphs.
Part of Hamidoune's main result of~\cite{hamidoune:hierarchical} is
obtained as a corollary. Additionally it is proven that
cycle-prefix graphs are optimally vertex connected.
Finally, for completeness, Mader's results~\cite{mader:minimale}
on the edge connectivity
of vertex transitive graphs and their proofs are given in
Section~\ref{section:edge connectivity}.

\section{Caley coset digraphs}
\label{section:caley coset}

Knowledge of basic group and graph theory is assumed (see for example
Herstein~\cite{grouptheory:text} and Tutte~\cite{graphtheory:text}).
All structures are assumed to be finite.
If $G$ is a group and $F$ is a union of left cosets $gH$ of a subgroup $H$
of $G$, then $F/H$ denotes the set of left cosets of $H$ in $F$.
We have $\bigcup(F/H)=F$. $\langle H_1,H_2,\ldots\rangle$
denotes the subgroup of $G$ generated by $H_1,H_2,\ldots$.

A digraph $\cG$ is {\em vertex transitive\/} iff the automorphisms of $\cG$
act transitively on the set of vertices of $G$.
The {\em transpose\/} of a digraph $G$, denoted by $G^*$,
is obtained by reversing all the edges of $G$.

Given a group $G$, a subgroup $H$ of $G$
and a set of {\em generators\/} $S\subseteq G\setminus H$,
the {\em Cayley coset digraph\/} $\cG(G,H,S)$ 
(or $\cG$ if $G$, $H$ and $S$ are clear from context) is obtained
as follows:
the set $V(\cG)$ of vertices of $\cG$ is given by the
set of left cosets $G/H=\{gH\suchthat g\in G\}$ of $H$ in $G$,
and the set $E(\cG)$ of edges
of $\cG$ consists of the ordered pairs $(gH, g'H)$
with $gHs\cap g'H\not=\emptyset$ for some
$s\in S$. The graph $\cG$
is vertex transitive. A transitive group of automorphisms of
$\cG$
acting on $V(\cG)$ consists of the maps
$\varphi_g:g'H\mapsto gg'H$.
Every vertex transitive digraph is a Cayley coset digraph, as is shown
in~\cite{sabidussi:thecomposition}.

$\cG$ is a {\em Cayley digraph\/} iff $H=\{e\}$, where
$e$ is the identity of $G$.
A {\em Cayley (coset) graph\/} is a symmetric Cayley (coset) digraph
(i.e. if $(x,y)\in E(\cG)$ then $(y,x)\in E(\cG)$).

For $s\in S$, let $$E_s=\{(g_1H,g_2H)\suchthat g_1Hs\cap
g_2H\not=\emptyset\}$$ be the set of edges
{\em induced\/} by $s$. The next lemma shows
that we can assume that $S$ consists
of distinct representatives of the double
cosets $HgH$ in $G$.

\begin{lemma}
\label{lemma:unique}%
The following are equivalent:
\begin{itemize}
\item[{\rm (i)}]The edge $(g_1H,g_2H)$ is in $E_s$.
\item[{\rm (ii)}]$g_1^{-1}g_2$ is in $HsH$.
\item[{\rm (iii)}]$g_1HsH\supseteq g_2H$.
\end{itemize}
\end{lemma}

\proof
We have $(g_1H,g_2H)\in E_s$ iff
$Hs\cap g_1^{-1}g_2H\not=\emptyset$ iff
$g_1^{-1}g_2\in HsH$ iff $g_2\in g_1HsH$ iff
$g_2H\subseteq g_1HsHH = g_1HsH$.
\qed

Lemma~\ref{lemma:unique} implies that $E_s\cap E_{s'}=\emptyset$  unless
$HsH=Hs'H$ in which case $E_s=E_{s'}$.

\noindpar{Assumption}From now
on we assume that the generators are representatives from
distinct double cosets of $H$.

A digraph $G$ is {\em strongly connected\/} iff for every $a,b\in V(G)$,
there is a path from $a$ to $b$. Since only strong connectivity
is considered in this note, the word ``strongly'' will be omitted.
Note that for vertex transitive digraphs, strong connectivity is
equivalent to weak connectivity, 
i.e. in a vertex transitive digraph, if there is a path from $a$ to
$b$, then there is one from $b$ to $a$.

\begin{lemma}
\label{lemma:gen-conn}%
$\cG$ is connected iff $\langle H,S\rangle=G$.
\end{lemma}

\proof
By vertex transitivity, it suffices to check
that for every $g\in G$, there is a path from $H$
to $gH$ in $\cG$. A path is a sequence
$H, g_1H, g_2H,\ldots, gH$ where
\[
g_1H\subseteq Hs_1H,\hspace{12pt} g_2H\subseteq g_1Hs_2H,\ldots
\]
with $s_i\in S$.
There is such a sequence ending at $gH$ iff every element
of $gH$ can be expressed as a product of elements of $H$ and generators.
The result follows.
\qed

Lemma~\ref{lemma:gen-conn}
implies that the components of $\cG$
are determined by the left cosets of $\langle H,S\rangle$, i.e.
each component is of the form $(g\langle H, S\rangle)/H$.

\noindpar{Assumption}We assume that
$H$ and $S$ generate $G$ or equivalently, that $\cG$ is connected.

By Lemma~\ref{lemma:gen-conn}, the set of neighbors of $H$ due to $s$
is given by $HsH/H$.
Thus, the contribution of $E_s$ to the degree $d$ of $\cG$
is determined by the index
$$d_s={\cardof{HsH/H}}.$$

\begin{lemma}
\label{lemma:degree}%
$d=\sum_{s\in S}d_s$.
\end{lemma}

\proof
The result follows from the discussion above and
from the assumption that
the generators come from distinct double cosets of $H$.
\qed

For $S'\subseteq S$, let $d_{S'}=\sum_{s\in S'}d_s$.

The numbers $d_s$ can be computed using the following elementary
result from group theory~\cite{grouptheory:text}:

\begin{lemma}
\label{lemma:cosetindex}%
$\cardof{HsH/H}=\cardof{H/(H\cap sHs^{-1})}$.
\end{lemma}

\section{Vertex connectivity}
\label{section:vertex connectivity}

\noindent{\bf Definitions.}
Let $G$ be a digraph. The vertex connectivity of $G$ is the smallest number
of vertices that need to be removed from $G$ so that the digraph induced
on the remaining vertices is not connected.
The vertex connectivity of $G$ is denoted by $\kappa(G)$
or simply $\kappa$ if the graph is clear from the context.

For $A\subseteq V(G)$, let $N_{G}(A)$ denote the set of neighbors
of $A$, where
$$N_{G}(A)= \{x\in V(G)\setminus A\suchthat
\exists y\in A\;\; \mbox{such that} \;\;(y,x)\in E(G)\}.$$
The subscript
is omitted if it is clear which graph is being considered.
$A$ is a {\em part\/} of $G$ iff
$V(G)\setminus(A\cup N(A))$ is non-empty.
$A$ is an {\em atom\/}  iff $A$ is a minimum size part of $G$
with the property that $\cardof{N(A)}=\kappa$.
Note that this definition 
differs slightly from the one given
by Hamidoune~\cite{hamidoune:sur}, who defines atoms as minimal size
parts $A$ of $G$
or $G^*$ satisfying $\cardof{N_G(A)}=\kappa$
or $\cardof{N_{G^*}(A)}=\kappa$, respectively.

The only digraphs without atoms are the complete digraphs
where every pair of vertices is an edge.
Assume that $\cG$ is not complete.
The next lemmas are used to show that
the atoms of $\cG$ partition $\cG$, provided
that they are small enough.

\begin{lemma}
\label{lemma:atomsdisj0}%
Let $G$ be a digraph on $n$ vertices of vertex connectivity $\kappa$.
Let $A$ be an atom of $G$ and $B$ a part of $G$ with
$\cardof{N(B)}=\kappa$.
If $V(G)\setminus(A\cup N(A)\cup B\cup N(B))$ is non-empty,
then either $A\cap B=\emptyset$ or $A\subseteq B$.
\end{lemma}

\proof
Suppose that $A\cap B$ and $A\setminus B$ are both non-empty.
We show that $A\cup B$ is a part with $\cardof{N(A\cup B)}<\kappa$
to derive a contradiction.
Since $V(G)\setminus( A\cup B\cup N(A\cup B))=
V(G)\setminus (A\cup N(A)\cup B\cup N(B))$, it follows that $A\cup B$
is a part. We have
$$N(A\cup B)=(N(B)\setminus A)\cup (N(A)\setminus (B\cup N(B)).$$
To obtain the contradiction,
it suffices to show that $\cardof{N(A)\setminus (B\cup N(B))}<\cardof{N(B)\cap A}$.
This follows from Lemma~\ref{lemma:simplec} which is proved next.
\qed

\begin{lemma}
\label{lemma:simplec}%
Let $G$ be a digraph with $A$ an atom and $B$ a part of $G$. If $A\cap B$
and $A\setminus B$ are both non-empty, then
$\cardof{N(A)\setminus (B\cup N(B))}<\cardof{N(B)\cap A}$.
\end{lemma}

\proof
Since $A$ is an atom, $\cardof{N(A\cap B)}>\cardof{N(A)}$.
Using the fact that
$N(A\cap B)$ is included in $(N(B)\cap A)\cup(N(A)\cap (B\cup N(B)))$,
we can deduce
\begin{eqnarray*}
\cardof{N(A)}&<&\cardof{(N(B)\cap A)\cup(N(A)\cap (B\cup N(B)))}\\
&=&\cardof{N(B)\cap A}+\cardof{N(A)\cap (B\cup N(B))}\\
&=&\cardof{N(B)\cap A}+\cardof{N(A)}-
\cardof{N(A)\setminus(B\cup N(B))}.
\end{eqnarray*}
This gives the result.
\qed

\begin{lemma}
\label{lemma:atomsdisj}%
Let $G$ be a digraph on $n$ vertices of vertex connectivity $\kappa$.
If the atoms of $G$ have size at most $(n-\kappa)/2$, then
any two distinct atoms of $G$
are disjoint.
\end{lemma}

\proof
By Lemma~\ref{lemma:atomsdisj0} it suffices to show that if
$A$ and $B$ are distinct intersecting atoms of $G$
of size at most $(n-\kappa)/2$, then
$V(G)\setminus(A\cup N(A)\cup B\cup N(B))\not=\emptyset$.
The set $A\cup N(A)\cup B\cup N(B)$ is the disjoint union of
$A\cup B$, $N(A)\setminus (B\cup N(B))$,
$N(A)\cap N(B)$ and $N(B)\setminus (A\cup N(A))$.
Thus
\begin{eqnarray*}
\cardof{A\cup N(A)\cup B\cup N(B)}&=&
\cardof{A\cup B}+\cardof{N(A)\cap N(B)}+\\
&&\hspace{10pt}
\cardof{N(A)\setminus (B\cup N(B))}+
\cardof{N(B)\setminus (A\cup N(A))}\\
&\leq& n-\kappa+\cardof{N(A)\cap N(B)}+\\&&\hspace{10pt}
\cardof{N(A)\setminus (B\cup N(B))}+
\cardof{N(B)\setminus (A\cup N(A))}.
\end{eqnarray*}
To bound the last sum, we use Lemma~\ref{lemma:simplec}
to obtain
\begin{eqnarray*}
2 \left[\cardof{N(A)\setminus(B\cup N(B))}\right.&+&
\left.\cardof{N(B)\setminus(A\cup N(A))}+
\cardof{N(A)\cap N(B)}\right]\\
&<&
\cardof{N(A)\cap B}+\cardof{N(B)\cap A}+\\&&\hspace{10pt}
\cardof{N(A)\setminus(B\cup N(B))}+
\cardof{N(B)\setminus(A\cup N(A))}+\\&&\hspace{10pt}
2\cardof{N(A)\cap N(B)}\\
&=&\cardof{N(A)}+\cardof{N(B)}\\
&=&2\kappa.
\end{eqnarray*}
It follows that $\cardof{A\cup N(A)\cup B\cup N(B)}< n-\kappa+\kappa=n$,
which gives the result.
\qed

\begin{lemma}
\label{lemma:atomsdisj1}%
Let $G$ be a digraph of vertex connectivity $\kappa$ which is not complete.
Then either $G$ or $G^*$ has an atom of size at most
$(n-\kappa)/2$.
\end{lemma}

\proof
Let $A$ be an atom of $G$. Then $ V(G)\setminus(A\cup N(A))$
is a part of $G^*$. Since $\cardof{N(A)}=\kappa$, the result follows.
\qed

\noindpar{Assumption}From now on we assume that $\cG$ has an atom of size
at most $(n-\kappa)/2$. To see that this assumption does not restrict
the generality of the results to be shown, it suffices to apply
Lemma~\ref{lemma:atomsdisj1} and note that the properties used to
prove the results are preserved if $\cG$ is replaced by $\cG^*$.  In
particular, note that the vertex connectivity of $\cG^*$ is the same
as the vertex connectivity of $\cG$. Furthermore, the generating set
$S^{-1}$ of $\cG^*$ has the same associated degrees and also consists
of distinct double coset representatives.

\begin{lemma}
\label{lemma:atoms-partition}%
With the given assumption,	
the atoms of $\cG$ partition the vertices
of $\cG$. The automorphisms of $\cG$ induce permutations of the atoms
and each atom is a vertex transitive induced subgraph of $\cG$.
\end{lemma}

\proof
The automorphic image of an atom is an atom.
By transitivity of $\cG$, the atoms cover $\cG$. Since distinct
atoms are disjoint, they partition the vertices of $\cG$.
This also implies that the automorphisms of $\cG$
induce permutations of the atoms. Since $\cG$ is transitive,
each atom is transitive.
\qed

Let $A_0$ be the atom containing $H$. Let $S_0$ be
the set of generators $s\in S$ such that
$s\in\bigcup A_0$. Let $S_1=S\setminus S_0$.	

\begin{lemma}
\label{lemma:atomssub}%
The subset $\bigcup A_0$ of $G$ is the
subgroup $\langle H, S_0\rangle$ generated
by $H$ and $S_0$. The edges induced on $A_0$ by $\cG$
are given by $E_{S_0}\cap (A_0\times A_0)$.
\end{lemma}

\proof
To see that $\bigcup A_0$ is a subgroup,
let $g\in\bigcup A_0$ and consider the automorphism
$\phi_g$. Since $\phi_g(H)\in A_0$, $\phi_g(A_0)=A_0$.
This implies that $\bigcup A_0$ is closed under multiplication
by $g$. That $\bigcup A_0$ is a subgroup follows
by arbitrariness of $g$.

If $g_1H$ and $g_2H$ are in $A_0$ and
there is an edge from $g_1H$ to $g_2H$ in $\cG$ induced by $s$,
then $g_2H\subseteq g_1HsH$.  In particular, there are
$h_1,h_2\in H$ such that $g_2=g_1hsh_2$. This gives
$h_1^{-1}g_1^{-1}g_2h_2^{-1} = s$.
Since $\bigcup A_0$ is
a subgroup of $G$, $s\in \bigcup A_0$,
so that $(g_1H,g_2H)\in E_{s_0}$.

Note that $A_0$ is connected, for
otherwise any part of $A_0$ with
outdegree zero in $A_0$ is a smaller part of $\cG$ with at most
$\kappa$ neighbors in $\cG$.
Lemma~\ref{lemma:gen-conn}
implies that $\cup A_0= \langle H, S_0\rangle$.
\qed

\begin{example}
Consider again the assumption that $\cG$ has an atom of size at most
$(n-\kappa)/2$.  There are Cayley graphs which do not satisfy this
assumption and where the atoms do not partition the set of vertices.
As an example consider the group $S_n$ of permutations on $n\geq 4$
vertices. Using cycle notation for permutations, let $a = (12)$ and $b
= (123\ldots n)$. Then $\langle a,b\rangle = S_n$. Let
$\cG=\cG(S_n,\{(1)\},\{a,b, ba\})$.  Then $S_n=\langle
a,b\rangle=\langle a, ba\rangle = \langle b, ba\rangle$. Therefore the
only candidates for atoms of $\cG$ or $\cG^*$ are
$H_1=\langle\rangle$, $H_2=\langle a \rangle$, $H_3 = \langle b
\rangle$ and $H_4 = \langle ba \rangle$.  Let $N_i=N(H_i)$ and
$N^*_i=N_{\cG^*}(H_i)$.  We have
\begin{eqnarray*}
N_1&=&\{a,b, ba\}\\
N^*_1&=&\{a,b^{-1},ab^{-1}\}\\
N_2&=&\{b,ba,ab,aba\}\\
N^*_2&=&\{ b^{-1}, ab^{-1}\}.\\
\end{eqnarray*}
Since $\cardof{H_3}=n$,
$\cardof{N_3}\geq n$
and $\cardof{N^*_3}\geq n$.
Since $ba$ is a cycle of length $n-1$, $\cardof{N_4}\geq n-1$ and
$\cardof{N^*_4}\geq n-1$.
Since $\cardof{N^*_2}=2$, $\kappa=2$. However,
none of the connected subgroups of $S_n$ have $2$ neighbors.
The atom containing the identity of $\cG^*$ is
given by $H_2$.
\end{example}

Lemma~\ref{lemma:atomssub} shows that to check the vertex connectivity
of $\cG$ it suffices to check the number of neighbors of
each subgroup generated by $H$ and a subset of $S$.
Since an atom is never the whole graph, the next result is immediate.

\begin{corollary}
\label{corollary:no-subgr}%
Suppose that for each $s\in S$, $H$ and $s$ generate $G$.
Then $\cG$ is optimally vertex connected.
\end{corollary}

Except when the degree of $\cG$ is $1$, the size
of atoms is strictly smaller than the degree $d$ of $\cG$.
This is a consequence of the next lemma.

\begin{lemma}
\label{lemma:atomsaresmall}%
$\cardof{N(A_0)}=\cardof{(\bigcup A_0)S_1H/H}\geq
\max(\cardof{A_0}, d_{S_1})$ and $\cardof{N(A_0)}$ 
is a multiple of $\cardof{A_0}$.
\end{lemma}

\proof
Let $s\in S_1$.
Then $HsH\cap(\bigcup A_0)=\emptyset$, for otherwise
$s$ induces an edge in $A_0$.
Since $\bigcup A_0$ is a subgroup, it follows
that $((\bigcup A_0)HsH)\cap(\bigcup A_0)=\emptyset$.
The set $(\bigcup A_0)HsH=(\bigcup A_0)sH$ is the union of
the neighbors of $A_0$ reachable by an edge in $E_s$.
Thus the union of the neighbors of $A_0$ is given by
$(\bigcup A_0)S_1H$. This implies the first identity. The inequality
is obtained by observing that $d_{S_1}=
\cardof{HS_1H/H}$. Since $(\bigcup A_0)S_1 H$ is a union
of right cosets of $(\bigcup A_0)$, $\cardof{N(A_0)}$ is a multiple
of $\cardof{A_0/H}=\cardof{A_0}$.
\qed

The following theorem generalizes Proposition~3.1 of~\cite{hamidoune:hierarchical}:

\begin{theorem}
\label{theorem:decomposition}%
Let $R_1$ and $R_2$ be a partition of $S$
such that the group $G'=\langle H,R_1\rangle$ does not
contain any members of $R_2$.
Let
$\cG'=\cG(G',H,R_1)$. Suppose that
if $r,s\in R_2$ and $G'rG'=G'sG'$, then
$\langle H,r\rangle=\langle H,s\rangle$.
Then $\kappa(\cG)\geq \min(\cardof{V(\cG')},\kappa(\cG')+d_{R_2})$.
\end{theorem}

\proof
The atom $A_0$ of $\cG$ satisfies one of the following cases:
\begin{itemize}
\item[1.]$\bigcup A_0\subset G'$,
\item[2.]$\bigcup A_0\supseteq G'$,
\item[3.]$\bigcup A_0$ and $G'$ are incomparable.
\end{itemize}

Consider case 1.
If $r\in R_2$, then $(\bigcup A_0)r H$ is disjoint from $G'$,
because otherwise $r\in G'$.
Thus the neighbors of $A_0$ due to $R_2$ are disjoint from $\cG'$ and
since $\cardof{(\bigcup A_0) rH/H}\geq\cardof{A_0}$,
there are at least $\max(\cardof{A_0},d_{R_2})$ many such neighbors.
It follows that
$$\cardof{N(A_0)}\geq\cardof{N_{\cG'}(A_0)}+\max(\cardof{A_0},d_{R_2}).$$
Either $\cardof{N_{\cG'}(A_0)}\geq\kappa(\cG')$ or
$A_0\cup N_{\cG'}(A_0)=V(\cG')$, and we are done.

Consider case 2.
Since the number of neighbors of $A_0$ is at least $\cardof{A_0}$,
trivially $\cardof{N_{\cG'}(A_0)}\geq \cardof{V(\cG')}$.

Consider case 3.
Let $A=A_0\cap V(\cG')$. 
Let $R_{11}=\{s\in R_1\suchthat s\not\in  \bigcup A\}$.
Let $R_{21}=\{s\in R_2\suchthat s\not\in \bigcup A_0\}$.
Let $R_{22}=\{s\in R_2\suchthat s\in\bigcup A_0\}$.
Observe that for $r\in R_{21}$ and $s\in R_{22}$,
$r$ and $s$ come from distinct double cosets of $G'$.
Otherwise, by assumption, $\langle H,r\rangle=\langle H,s\rangle$, which
would imply that either both $r$ and $s$ are in  $\bigcup A_0$,
or both $r$ and $s$ are not in $\bigcup A_0$.
Thus $G'R_{21}G'$ and $G'R_{22}G'$ are disjoint.
We have $(\bigcup A)R_{22}H\subseteq \bigcup (A_0\setminus V(\cG'))$,
so that the set $N_1=(\bigcup A)R_{22}HR_{11}H/H$
consists of neighbors of $A_0$ not in $\cG'$.
The set $N_2=(\bigcup A)R_{21}H/H$ also
consists of
neighbors of $A_0$ not in $\cG'$. Since
$(\bigcup A)R_{22}HR_{11}H\subseteq G'R_{22}G'$
and $(\bigcup A)R_{21}H\subseteq G'R_{21}G'$,
$N_1$ and $N_2$ are disjoint. It follows
that
$$\cardof{N_{\cG}(A_0)}\geq
\cardof{N_{\cG'}(A)}+\cardof{N_1}+\cardof{N_2}.$$
If $N_{\cG'}(A)\cup A=V(\cG')$, then  using
$\cardof{N_1}\geq \cardof{A}$ we get
$\cardof{N_{\cG}(A_0)}\geq\cardof{V(\cG')}$.
If not, then using $\cardof{N_1}\geq d_{R_{22}}$
and $\cardof{N_2}\geq d_{R_{21}}$ gives
$$\cardof{N_{\cG}(A_0)}\geq\cardof{N_{\cG'}(A)}+d_{R_{22}}+d_{R_{21}}\geq
\kappa(\cG')+d_{R_2},$$ as desired.
\qed

\begin{corollary}
\label{corollary:1}%
Let $\{S_1,\ldots, S_k\}$ be a partition of $S$.
Define $G_i=\langle H,S_1,\ldots, S_i\rangle$ and
$d_i=d_{S_1}+\ldots+d_{S_i}$.
Suppose that the following hold:
\begin{itemize}
\item[1.]The $G_i$ are distinct.
\item[2.]If $r,s\in S_{i+1}$ and $G_irG_i=G_isG_i$,
then $\langle H,r\rangle=\langle H,s\rangle$.
\item[3.]$\kappa(\cG(G_1,H,S_1))=d_1$.
\item[4.]$\cardof{G_i/H}\geq d_{i+1}$.
\end{itemize}
Then $\kappa(\cG)=d_S$.
\end{corollary}

\proof
It suffices to apply Theorem~\ref{theorem:decomposition} to
each step in the tower of Cayley coset digraphs
$\cG(G_i,H,S_1\cup\ldots\cup S_i)$.
\qed

The next corollary
shows how to replace the restriction on the size of the groups $G_i$ by
conditions on the degrees induced by the partition of the generators.

\begin{corollary}
\label{corollary:1.1}%
Let $\{S_1,\ldots S_k\}$ be a partition of $S$. Define $G_i=\langle
H,S_1,\ldots S_i\rangle$ and $d_i=d_{S_1}+\ldots+d_{S_i}$.
Suppose that the following hold:
\begin{itemize}
\item[1.]The $G_i$ are distinct.
\item[2.]If $r,s\in S_{i+1}$ and $G_irG_i=G_isG_i$,
then $\langle H,r\rangle=\langle H,s\rangle$.
\item[3.]$\kappa(\cG(G_1,H,S_1))=d_1$.
\item[4.]$\cardof{G_1/H}\geq d_2$.
\item[5.]For all $i$, $d_{S_{i+1}}\leq d_i$.
\end{itemize}
Then $\kappa(\cG)=d_S$.
\end{corollary}

\proof
It suffices to show by induction that $\cardof{G_i/H}\geq d_{i+1}$
and apply Corollary~\ref{corollary:1}. 
Since $G_{i-1}$ is a proper subgroup of $G_{i}$,
$\cardof{G_i/H}\geq 2\cardof{G_{i-1}/H}\geq 2 d_i$.
Since $d_{S_{i+1}}\leq d_i$, the result follows.
\qed

\section{Applications}
\label{section:applications}

\noindent{\bf Hierarchical Cayley coset digraphs.}

\begin{definition}
$\cG$ is a {\em quasi-minimal\/} or {\em hierarchical\/}
Cayley coset digraph iff there
is an ordering $\{s_1,\ldots,s_k\}$ of the generators of $\cG$
such that the subgroups $\langle H,s_1,s_2,\ldots,s_i\rangle$
are distinct.
$\cG$ is {\em minimal\/} iff for no $S'\subset S$, $\langle H, S'\rangle = G$.
\end{definition}

Corollary~\ref{corollary:1.1} can be simplified for hierarchical Cayley coset digraphs.

\begin{theorem}
\label{theorem:hierarchical-gen}%
Let $\cG$ be hierarchical with generators
ordered by $\{s_1,\ldots,s_k\}$. Let $d_i=d_{s_1}+\ldots+d_{s_i}$
and $G_i=\langle H,s_1,\ldots, s_i\rangle$.
Suppose that for each $i$, $d_{s_{i+1}}\leq d_i$
and  $\cardof{G_1/H}\geq d_2$.
Then $\kappa(\cG)=d(\cG)$.
\end{theorem}

\proof
Condition 1 of Corollary~\ref{corollary:1.1} is satisfied by the
definition of hierarchical Cayley coset digraphs.  For the partition
of $S$ into singletons, condition 2 is trivially satisfied.  A Cayley
coset digraph generated by a single generator is optimally vertex connected
(Corollary~\ref{corollary:no-subgr}), so that condition 3 is
satisfied.  Conditions 4 and 5 are satisfied by assumption.  \qed

\begin{corollary}
\label{corollary:hier1}%
The assumption that $\cardof{G_1/H}\geq d_2$ in the statement
of Theorem~\ref{theorem:hierarchical-gen} can be replaced by
the assumption that
$Hs_1^{-1}H\not=Hs_1H$.
\end{corollary}

\proof
The assumption that $Hs_1^{-1}H\not=Hs_1H$ is equivalent to the
assumption that there are no cycles of length two in $\cG_1$.
This implies $\cardof{V(\cG_1)}\geq 2 d_1+1>d_2$.
\qed

The fact that hierarchical Cayley digraphs are optimally vertex connected 
can now be easily shown.

\begin{corollary}
\label{corollary:hierarchical}%
{\rm (Baumslag~\cite{baumslag:fault}, Hamidoune~\cite{hamidoune:hierarchical})}
Hierarchical Cayley digraphs are optimally vertex connected.
\end{corollary}

\proof
Let $\cG$ be a hierarchical Cayley digraph with generators ordered
by $\{s_1,\ldots, s_k\}$.
Since $\cardof{\langle s_1\rangle}\geq 2$ and $d_{s_i}=1$
for each $i$,
the result is immediate by Theorem~\ref{theorem:hierarchical-gen}.
\qed

Part of the more general result in Hamidoune~\cite{hamidoune:hierarchical} also
follows and was obtained by Hamidoune using a restricted version
of Theorem~\ref{theorem:decomposition}.

\begin{theorem}
\label{theorem:hierarchical-gen-c}%
{\rm (Hamidoune~\cite{hamidoune:hierarchical})}
Let $\cG=\cG(G,\{e\}, S\cup S')$, where $S'\subseteq S^{-1}$ and the
elements of $S'$ have order at least three. Assume that
$\cG(G,\{e\}, S)$ is hierarchical with the elements of $S$ ordered
by $S=\{s_1,\ldots, s_k\}$. 
If $\cardof{\langle s_1,s_2\rangle}\not=4$ then $\cG$ is optimally
vertex connected.
\end{theorem}

\proof
Partition $S\cup S'$ into the sets defined by: 
$S_1= \{s_1,s_2\}\cup (\{s_1^{-1},s_2^{-1}\}\cap S')$ and for $i>1$,
$S_i=\{s_{i+1}\}$ if $s_{i+1}^{-1}\not\in S'$ and
$S_i=\{s_{i+1},s_{i+1}^{-1}\}$ otherwise.
Let $d_i$ and $G_i$ be defined as in Corollary~\ref{corollary:1.1}.
We have $2\leq d_1\leq 4$, and
for each $i\geq 1$, $d_{S_{i+1}}\leq 2$.
To apply Corollary~\ref{corollary:1.1} it suffices to check
that $\cG_1=\cG(G_1,\{e\},S_1)$ is optimally vertex
connected and $\cardof{G_1}\geq d_2$.

If $\cG_1$ is not optimally vertex connected, then the atom $A$ 
of $\cG_1$ is given by
either $\langle s_1\rangle$ or $\langle s_2\rangle$.  Let $s_i\not\in
A$ and $s_j\in A$ ($\{i,j\}=\{1,2\}$).  We can assume that $2\leq
\cardof{A}\leq 3$, for otherwise $\cardof{N(A)}\geq\cardof{A}\geq
4\geq d_1$.  Similarly, since $\cardof{A}$ divides $\cardof{N(A)}$, we
can assume that $\cardof{N(A)}=\cardof{A}$.  Suppose that $s_i^{-1}\in
S'$ so that $s_i\not=s_i^{-1}$.  Then the neighbors of $A$ are given
by expressions of the form $s_j^{l}s_i$ and $s_j^{l}s_i^{-1}$.  Since
$\cardof{N(A)}=\cardof{A}$, for some $l\not=l'$,
$s_j^ls_i=s_j^{l'}s_i^{-1}$.  This implies that $s_i^2=s_j^{l'-l}$ so
that $A\cup N(A)=G_1$, contradicting the assumption that $A$ is an
atom. Thus $s_i^{-1}\not\in S'$. This implies that $d_1\leq 3$ so that
$\cardof{A}=2$.  But this implies that $s_j^{-1}\not\in S'$, whence
$d_1=2$ contradicting the assumption that $A$ is a nontrivial atom.

By assumption, $G_1$ has a proper nontrivial subgroup and
$\cardof{G_1}\not=4$. This implies that $G_1\geq 6\geq d_2$.
\qed

Hamidoune~\cite{hamidoune:hierarchical} continues the analysis of the
proof of Theorem~\ref{theorem:hierarchical-gen-c} to show that if
$\cG$ is as in the statement of this theorem and $\cG$ is not
optimally vertex connected, then $k\geq 3$, $\kappa(\cG)=\cardof{S\cup
S'}-1$, $s_i^2=s_1$ for $i>1$, $s_1^2=1$ and
$S'=(S\setminus\{s_1\})^{-1}$.

\noindent{\bf Cycle-prefix graphs.}
The cycle-prefix graphs (CP-graphs) are Cayley coset digraphs defined
on the group $S_n$ of permutations of $[n]=\{1,\ldots,n\}$. They
were proposed as interconnection networks with good
degree and diameter properties in~\cite{chen:basicCP}. If $\pi$
is a permutation which maps $i$ to $\pi_i$, then we write
$\pi=\pi_1\pi_2\ldots\pi_n$. Application of permutations is on the
right, so that $i\pi=\pi_i$. Composition is defined by
$i(\pi\sigma)=(i\pi)\sigma$.  The cycle-prefix permutations
$\gamma(k)$ with $2\leq k\leq n$ are the permutations which cyclically
permute $\{1,\ldots,k\}$ to the right and leave other numbers fixed.
Thus $\gamma(k)=k12\ldots(k-1)(k+1)\ldots(n-1)n$.  Let $H_k$ be the
subgroup of $S_n$ consisting of the permutations $\pi$ with $\pi_i=i$
for $i\leq n-k$.  For $1\leq k\leq n-1$, let
$\mfn{CP}(n,k)=\cG(S_n,H_k,\{\gamma(2),\ldots,\gamma(n-k+1)\})$.  Then
$\mfn{CP}(n,k)$ is a Cayley coset digraph of degree $n-1$.  The
degrees induced by the generators in $\mfn{CP}(n,k)$ are given by
$d_{\gamma(i)} = 1 $ for $2\leq i \leq n-k$ and $d_{\gamma(n-k+1)} =
k$.

\begin{theorem}
\label{theorem:cp-connected}
$\mfn{CP}(n,k)$ is optimally vertex connected.
\end{theorem}

\proof
The digraphs $\mfn{CP}(n,k)$ are hierarchical with generators ordered
by $\{\gamma(2),\ldots,\gamma(n-k+1)\}$.  If $k=n-1$, then
$\mfn{CP}(n,k)$ is the complete digraph on $n$ vertices and we are
done.  If $k=1$, then $\mfn{CP}(n,k)$ is a hierarchical Cayley digraph
and is optimally vertex connected by Corollary~\ref{corollary:hierarchical}.
Assume that $1<k<n-1$.  Let $G'=\langle H,
\gamma(2),\ldots,\gamma(n-k)\rangle$ and
$\cG'=\cG(G',H,\{\gamma(2),\ldots,\gamma(n-k)\})$.  Then $\cG'$ is
isomorphic to $\mfn{CP}(n-k,1)$ and $\cardof{\cG'}=(n-k)!$.  This
follows from the fact that $\gamma(j)$ is in the normalizer of $H$ for
$j\leq n-k$.  Thus $\cG'$ is optimally vertex connected.  If $(n-k)!\geq
n-1$, Theorem~\ref{theorem:decomposition} can be applied to show that
$\mfn{CP}(n,k)$ is optimally vertex connected.  Assume that
$(n-k)!<n-1$. Then $k\geq n/2$, because $\lceil (n+1)/2\rceil!\geq
n-1$ for all $n\geq 1$.  Following the proof of
Theorem~\ref{theorem:decomposition}, suppose that $A$ is a nontrivial
atom of $\mfn{CP}(n,k)$ with $\cardof{A}>1$ and $\cardof{N(A)}< n-1$.
Note that $\langle H,\gamma(n-k+1)\rangle = S_n$ which implies that
the edges due to $\gamma(n-k+1)$ are not in $A$ and $\bigcup
A\subseteq G'$.  For any subgroup $F$ of $G'$ which contains $H$, the
number of neighbors of $F/H$ due to $\gamma(n-k+1)$ is $\cardof{F/H}
k$.  Proof: $F\gamma(n-k+1)H$ is given by the disjoint union
$\bigcup_{l=0}^{k-1}F_l$,
where $F_l$ is the set of permutations defined by
\[
F_l=\{ \pi\suchthat \mbox{$\pi_1=n-l$ and $\exists\sigma\in F$ such that
for $2\leq i\leq n-k$, $\pi_i=\sigma_{i-1}$}\}.
\]
Note that for each member $\sigma$ of $F$, $\sigma_{n-k}$ is determined
by the $\sigma_1\ldots\sigma_{n-k-1}$. Hence $\cardof{F_l}=\cardof{F}$
which gives $\cardof{F\gamma(n-k+1)H/H} = \cardof{F/H}k$.
As a result, $\cardof{A}<n/k$ and
since $\cardof{A}>1$ this implies that $k<n/2$, contrary to assumption.
\qed

\section{Edge connectivity}
\label{section:edge connectivity}

For completeness we include the result that vertex transitive
digraphs are optimally edge 
connected. This result and its proof (for the undirected case)
are due to Mader~\cite{mader:minimale}.

\noindent{\bf Definition.}
The edge connectivity of a digraph $G$ is the smallest number
of edges that need to be removed so that the resulting digraph
is not connected. The edge connectivity of $G$ is denoted
by $\lambda(G)$.

\begin{theorem}
\label{theorem:edgec}%
Every Cayley coset digraph has edge connectivity equal to its degree.
\end{theorem}

\proof
Let $\lambda$ be the edge connectivity of $\cG$.
We start by showing that there is a notion of atom applicable to
edge connectivity. An {\em e-atom\/} of $\cG$
is a minimal subset of the vertices of $\cG$ with exactly $\lambda$ outgoing edges.
Let $N_e(A)$ denote the set of edges leaving $A$.
Let $E(A)$ denote the set of edges included in $A$.

\begin{lemma}
\label{lemma:eatom}%
If $A$ is an $e-atom$ and $B\subset V(\cG)$ with
$\cardof{N_e(B)}=\lambda$,
then $A\subseteq B$ or $A\cap B=\emptyset$ or
$A\cup B = V(G)$.
\end{lemma}

\proof
Suppose that $A\cap B\not=\emptyset$ and $A\cup B \not= V(G)$.
\begin{eqnarray*}
N_e(A\cap B)&=& (N_e(A)\cap E(B))\cup (N_e(A)\cap N_e(B))\cup(N_e(B)\cap E(A)),
\\
N_e(A\cup B)&=&( N_e(A)\setminus (E(B)\cup N_e(B)))\cup N_e(B)\setminus E(A),
\end{eqnarray*}
where the unions are disjoint.
This gives
\begin{eqnarray*}
\cardof{N_e(A\cup B)} &=& \cardof{N_e(A)\setminus (E(B)\cup N_e(B))}+
  \cardof{N_e(B)\setminus E(A)}\\
 &=&\cardof{N_e(A)\setminus (E(B)\cup N_e(B))} -
  \cardof{N_e(B)\cap E(A)} + \cardof{N_e(B)}\\
 &=& \cardof{N_e(A)\setminus (E(B)\cup N_e(B))} -
  \cardof{N_e(B)\cap E(A)} + \lambda.
\end{eqnarray*}
Since $\cardof{N_e(A\cup B)}\geq\lambda$, 
$$\cardof{N_e(A)\setminus (E(B)\cup
N_e(B))}\geq \cardof{N_e(B)\cap E(A)}.$$
Hence
$$
\cardof{N_e(A\cap B)}\leq \cardof{N_e(A)\cap E(B)}+
\cardof{N_e(A)\cap N_e(B)}+\cardof{N_e(A)\setminus (E(B)\cup N_e(B))}.
$$
Since the right-hand side of this expression is $\cardof{N_e(A)}=\lambda$,
minimality of $A$ requires that $A\subseteq B$.
\qed

Lemma~\ref{lemma:eatom} implies that the observations about
atoms of $\cG$ apply to e-atoms of $\cG$. In particular,
provided that the size of an e-atom is at most $n/2$,
distinct e-atoms are disjoint, so that they form blocks under
the group of automorphisms of $\cG$. 
Again, either $G$ or $G^*$ satisfies the size condition, so without loss
of generality, assume that an atom has at most
$n/2$ elements.

Let $A_0$ be the e-atom which contains $H$. Then $\bigcup A_0$
is a subgroup of $G$. We can partition the generators as
before into the set $S_0$ of generators in $\bigcup A_0$
and $S_1$ of generators outside $\bigcup A_0$.
In this case the analysis is simple:
$A_0$ contains at least $d_{S_0}+1$ members each of which
has at least $d_{S_1}$ edges going outside of $A_0$.
Thus $\cardof{N_e(A_0)}\geq d_{S_1}(d_{S_0}+1)\geq d_{S_0}+d_{S_1}= d_S$.
Thus $\cardof{A_0}=1$ and
$\lambda=\cardof{N_e(A_0)}=d_S$.
\qed

\noindpar{Acknowledgements}
This work would not have been possible without the helpful discussions
with Vance Faber and Bill Chen.

\end{document}